\documentclass[12pt]{article}
\usepackage{amsmath} 
\usepackage{amsbsy} 
\usepackage{amssymb,amsthm,psfig}

\newtheorem{theorem}{Theorem}[section] 
 
\newtheorem{lemma}{Lemma}[section] 

 \newtheorem{definition}{Definition}[section] 
\newtheorem{remark}{Remark}[section] 
  
 \numberwithin{equation}{section} 
 \newcommand{\norm}[1]{\left\Vert#1\right\Vert}

\title{Unconditional Stability for Numerical Scheme Combining Implicit Timestepping for Local Effects and Explicit Timestepping 
for Nonlocal Effects}
\author{Mihai Anitescu\\ 
        anitescu@mcs.anl.gov\\ 
        Argonne National Laboratory\\ 
        Argonne, IL 60439, U.S.A\\ 
        William J. Layton \\
        and Faranak Pahlevani \\
        wjl@pitt.edu, fap4@pitt.edu\\ 
        Department of Mathematics\\ 
        University of Pittsburgh\\ 
        Pittsburgh, PA 15260, U.S.A} 
 
\begin{document} 
\maketitle 

\vspace{-6in} {\tiny Preprint ANL/MCS-P1093-0903} 
\vspace{6.2in} 

\begin{abstract} 
A combination of implicit and explicit timestepping is analyzed  
for a system of ODEs motivated by ones arising from spatial discretizations of 
evolutionary partial differential equations. Loosely speaking, the 
method we consider is implicit in local and stabilizing terms in 
the underlying PDE and explicit in nonlocal and unstabilizing 
terms. Unconditional stability and convergence of the 
numerical scheme are proven by the energy method and by algebraic techniques.  
This stability result is surprising because usually when different methods are  
combined, the stability properties of the least stable method plays a determining  
role in the combination. 
\end{abstract} 
 
\section{Introduction} 
This report considers timestepping methods for systems of ordinary 
differential equations of the form 
\begin{equation} \label{eq:mainEquation} 
 u'(t) + Au(t)+B(u)u(t)-Cu(t)=f(t), 
\end{equation} 
in which $A$, $B(u)$, and $C$ are ${n}\times{n}$  
matrices, $u(t)$ and $f(t)$ are $n$-vectors, and  
\begin{equation}\label{eq:condEquation} 
 A=A^{T} \succ 0,\, B(u)=-B(u)^{*},\, C=C^{T}\succeq 0 \mbox{~and~} A-C \succeq 0. 
\end{equation} 
Here $\succ$ and $\succeq$ denote, respectively, the positive definite and  
the positive semidefinite ordering.  
The key properties motivating our work are that $A$ is sparse and 
that although $C$ is not sparse, the action of $C$ on a vector is 
inexpensive to calculate. This structure is motivated by 
multiscale discretizations of turbulence but can also arise from 
closed-loop control problems and ensemble calculations. Given this  
structure of (\ref{eq:mainEquation}), the simplest scheme that is computationally feasible is 
\underline{explicit} in the global, unstable part of (\ref{eq:mainEquation}), that is, 
$Cu$. Accordingly, we consider 
\begin{equation} \label{eq:schemeEquation} 
\frac{u_{n+1}-u_n}{k}+A{u_{n+1}}+B({u_n}){ 
u_{n+1}}-C{u_n}=f_{n+1}, k=\Delta t, 
\end{equation} 
where $u_n$ is the approximation to $u(t=nk)$.
Usually when methods are combined, the stability properties of the 
explicit method play a determining role in the overall method. In 
Theorems \ref{t:unconditional} and \ref{t:bouHom}, we prove the surprising result that 
(\ref{eq:schemeEquation}) is \underline {unconditionally stable}. This result is outside the realm of  
root condition stability analysis for uncoupled scalar problems.
 
In Section 2, unconditional stability and convergence of (\ref{eq:schemeEquation}) are proven. 
We give two stability proofs. The first is algebraic. Since the constants 
depend on the dimension of the system, we also give an energy 
proof of stability (with uniform constants) that is potentially 
extensible to discretized PDEs. Section 3 presents numerical tests 
illustrating the theory. First, we briefly summarize some 
motivating problems leading to (\ref{eq:mainEquation}). 
 
The basic model of the turbulent dispersion is that it is 
dissipative in the mean (see ~\cite{KP80},~\cite{MP93},~\cite{IL98}). 
A more accurate formulation is that its 
dissipative effects are focused on the smallest resolved scales
(see ~\cite{HMJ00}). 
This physical idea has led to algorithms for numerical 
stabilization of transport-dominated phenomena based on eddy 
diffusivity acting only on the smallest resolved scales (e.g., 
~\cite{KA02},~\cite{HU00},~\cite{TA89},~\cite{KL02},~\cite{G99a},
~\cite{G99b},~\cite{HMJ00},~\cite{L00},~\cite{L02}). The natural
realization of this idea for spatial discretizations of   
convection diffusion equations is diffusive stabilization on all 
scales and then 
antidiffusing on the large scales. This leads to the system of ODEs 
\begin{equation} \label{eq:tensor}
\dot u_{ij}(t)+b \cdot \nabla^h u_{ij} - 
(\epsilon_0(h)+\epsilon) \Delta^h u_{ij} + \epsilon_0(h) 
P_H(\Delta^h P_H(u_{ij}))= 
f_{ij}, 
\end{equation} 
where standard notation is used: $\Delta^h$ is the discrete 
Laplacian, $\epsilon_0(h)$ is the artificial viscosity parameters and 
$P_H$ denotes a projection onto a coarser mesh; see Section 
3 for details. The system (\ref{eq:tensor}) fits exactly the form (\ref{eq:mainEquation}), 
(\ref{eq:condEquation}), where $C$ is provided   as the matrix arising from  
$\epsilon_0(h)$ term. We shall also test one algorithm as a 
perturbation of the method (\ref{eq:tensor}) in which the projection is replaced 
by a nearest averaging $\overline{\Delta^h \overline{{u}_{ij}}}$. In both cases,
the projection or averaging operator 
accounts for the nonlocal character (i.e., the large bandwidth) 
of $C$. On the other hand, averaging and 
projection are both embarrassingly parallel 
operators whose action on a given vector is cheap to perform.\\ 
\begin{remark} (1) A second main application is discretization of 
turbulent flow problems which, although nonlinear and constrained, 
have a similar structure to the above (simple) 
linear convection diffusion problem.\\ 
(2) A known method of stabilizing the timestepping and the 
associated linear system (but not the spacial discretization ) 
corresponds to (\ref{eq:mainEquation}) without the averaging: 
\begin{equation} 
\frac{u_{n+1}-u_n}{k}+b \cdot \nabla^h u_{n+1} - 
(\epsilon_0(h)+\epsilon) \Delta^h u_{n+1} + \epsilon_0(h) \Delta^h 
u_{n}= f_{n+1}. 
\end{equation} 
Each time step requires the inversion of the matrix corresponding 
the operator $- (\epsilon_0(h)+\epsilon) \Delta^h+b \cdot 
\nabla^h+k^{-1}I$, which, for $\epsilon_0$ suitably chosen, is an 
$M$-matrix. Our analysis applies to this method as well. 
\end{remark} 
 
\section{The Stability Analysis} 
For our analysis, we assume that $B(u)$ is in $C^1(\Re^n)$ and 
$f(t)$ is in $C^1([0,\infty))$. For any $T>0$, we denote by 
\[ F_T=\max_{t \in [0,T]} \norm{f(t)}_2. \]

\begin{lemma} \label{l:boundedSolution}
The system of ODEs (\ref{eq:mainEquation}) under the 
condition (\ref{eq:condEquation}) with initial condition 
$u(0)=u_0$ has a  unique solution on $[0,T]$, for any $T>0$. 
\end{lemma}

{\bf Proof} Since (\ref{eq:mainEquation}) can be written as 
$\dot{u}=\psi(t,u)$ with $\psi$ being of class 
$C^0$ in $t$ and $C^1$ in $u$, local existence and 
uniqueness follows from the classical theory of ODEs \cite[Theorem V.8]{BR62}. 

We now show that the solution does not experience blow-up and can 
be extended everywhere. We multiply through (\ref{eq:mainEquation}) by 
$(u(t)^T)$ and we use (\ref{eq:condEquation}) to obtain 
that 
\[ u(t)^T u'(t) \leq -u(t)^T (A-C) u(t) +  u(t)^T f(t) \leq u(t)^T f(t).\]
Using Cauchy-Schwarz, we obtain that 
\[ \frac{d}{dt} \norm{u (t)}_2^2 \leq \norm{u(t)}_2^2 + F_T^2. \]
In turn, this implies that, 
\[ \norm{u(t)}_2^2 \leq \norm{u(0)}_2^2 e^t + F_T^2 \left(e^t-1\right) \] 
for any $t$ in an interval containing $0$ where $u(t)$ is defined. Since 
$u(t)$ does not experience blow-up in finite time, it can 
be extended uniquely over all of $[0,T]$. $\Box$

Note that from (\ref{eq:mainEquation}) and our assumtion that
$f(t)$ is of class $C^1([0,\infty])$, we get that $u(t)$ is of class
$C^2([0,\infty])$. The fact that $u''(t)$ is continuous will be used 
in determining a bound for the truncation error. 

Consider the system of ODEs (\ref{eq:mainEquation}) under the 
condition (\ref{eq:condEquation}) and discretized
by (\ref{eq:schemeEquation}). 

First, we note that each step of (\ref{eq:schemeEquation}) requires
the inversion of $I+kA+kB_n$.   
\begin{lemma} 
Under (\ref{eq:condEquation}) the ${n}\times{n}$ 
matrix $I+kA+kB_n$ has a positive definite symmetric part and 
is invertible. 
\end{lemma}  
\noindent 
{\bf Proof:} Let $x$ be any nonzero vector in $\Re^n$. Then 
\[ \begin{array}{rcl} x^T(I+kA+kB_n)x & = & x^Tx+kx^TAx+kx^TB_nx \\
 & = & {\norm{x}_2}^2+kx^TAx>0.~~~
\Box{}{} \end{array} \]
Since $A$, $B_n=B(u_n)$ and $C$ do not commute, the stability of the 
numerical scheme cannot be analyzed by reduction to eigenvalues.
Therefore, we formulate an energy norm  
that is not increased at each time step, that is, $\norm{u_{n+1}}_E \le  
\norm{u_n}_E$. 
\begin{definition} 
The energy norm of (\ref{eq:schemeEquation}), $\norm{.}_E$, is given by  
\begin{equation} \label{def:ENorm}
{\norm{u}_E}^2={u}^Tu +k{u}^TCu,
\end{equation}
for some $u \in \Re^n$, and its associated  
inner product is $<u,v>_E=(N_kv)^T(N_ku)$, with 
$N_k=(I+kC)^{\frac{1}{2}}$, for some $u,v \in \Re^n$. 
\end{definition} 

It can be seem immediately that the energy norm and the 2-norm satisfy the 
following inequality:

\[ \sqrt{1 + k \lambda_{min}(C)} \leq \norm{u}_E \leq  \sqrt{1 + k \lambda_{max}(C)}, \]
where $\lambda_{min}(C)$ and $\lambda_{max}(C)$ are, respectively, the 
smallest and the largest eigenvalue of $C$. From this inequality  and the positive
semidefiniteness of $C$, we get that the induced matrix norms satisfy
\[ \norm{A}_E \leq \norm{A}_2 \sqrt{1 + k \lambda_{max}(C)}. \]

\begin{theorem} \label{t:unconditional} 
Let $u_n$ satisfy (\ref{eq:schemeEquation}) with ${f(.)}\equiv{0}$,
under the condition (\ref{eq:condEquation}) on the coefficients. Then,
\[\norm{u_{n+1}}_E \le \norm{u_n}_E.\] 
\end{theorem} 
\noindent 
{\bf Proof:} Multiplying with $u_{n+1}^T$ through the equation in (\ref{eq:schemeEquation}), we obtain  
\[u_{n+1}^T\frac{u_{n+1}-u_n}{k}+u_{n+1}^TA{u_{n+1}}+u_{n+1}^TB_n 
u_{n+1}=u_{n+1}^TCu_n\]
Since $B_n$ is skew symmetric, $u_{n+1}^TB_nu_{n+1}=0$. Therefore 
\begin{equation} 
u_{n+1}^T\frac{u_{n+1}-u_n}{k}+u_{n+1}^TA{u_{n+1}}=u_{n+1}^TCu_n. 
\end{equation} 
This is equivalent to 
\begin{equation} 
u_{n+1}^Tu_{n+1}+ku_{n+1}^TAu_{n+1}=ku_{n+1}^TCu_n+u_{n+1}^Tu_n .
\end{equation} 
Since $A \succeq C$, we have that
\begin{equation} \label{eq:Cineq}
u_{n+1}^Tu_{n+1}+ku_{n+1}^TCu_{n+1} \le u_{n+1}^Tu_n 
+ku_{n+1}^TCu_n. 
\end{equation} 
Define $w=(u_{n+1},k^{1/2}C^{1/2}u_{n+1})^T,
v=(u_n,k^{1/2}C^{1/2}u_n)^T$. Then (\ref{eq:Cineq})) can be written as $w^Tw 
\le w^Tv$. Applying the Cauchy-Schwarz inequality, we get $\norm{w}_2  
\le \norm{v}_2$. Hence, 
\begin{equation} 
u_{n+1}^Tu_{n+1}+ku_{n+1}^TCu_{n+1} \le u_n^Tu_n +ku_n^TCu_n 
\end{equation} 
or
\[\norm{u_{n+1}}_E \le \norm{u_n}_E. \Box{}{}\]  

The conclusion of the preceding theorem is that when 
(\ref{eq:mainEquation}) is homogeneous, $f \equiv 0$, we obtain that 
$\norm{u_n}_E \neq \norm{u_0}_e$, $\forall n$, independent of $T$. 
This means that our method is, indeed, \underline{unconditionally stable}. 

Consider (\ref{eq:schemeEquation}) with ${f}\equiv{0}$, rewritten as
\begin{equation} \label{eq:homEquation}
(I+kA+kB_n)u_{n+1}=(I+kC)u_n, B_n=B(u_n).
\end{equation}  
Equation (\ref{eq:homEquation}) yields 
\[
u_{n+1}=(I+kA+kB_n)^{-1}(I+kC)u_n, 
\]
which, in turn, implies that 
\[
(I+kC)^\frac{1}{2}u_{n+1}=(I+kC)^\frac{1}{2} (I+kA+kB_n)^{-1}(I+kC)^\frac{1}{2} (I+kC)^\frac{1}{2} u_n. 
\]
Therefore, from the definition of $\norm{\cdot}_E$, a sufficient 
condition to prove the unconditional stability result is to prove that
\[ 
\norm{(I+kC)^\frac{1}{2} (I+kA+kB_n)^{-1}(I+kC)^\frac{1}{2}}_2 \leq 1. 
\]
From \ref{eq:condEquation}, this can be done by using the following Lemma. 

\begin{lemma} \label{l:mainLemma} 
Let $D_1={D_1}^T \succ 0$ and $D_2={D_2}^T \succ 0$ be ${n}\times{n}$ matrices such that $D_1-D_2 \succ 0$. Let 
$D_4=D_2^{\frac{1}{2}}$ and be symmetric. 
If $D_3$ is an ${n}\times{n}$ skew-symmetric matrix, then 
\begin{equation} 
\parallel D_4(D_1+D_3)^{-1}D_4 \parallel _2 \le 1. 
\end{equation} 
\end{lemma} 
\noindent 
{\bf Proof:} Let $F=D_4(D_1+D_3)^{-1}D_4$. It is straightforward that\\ 
$F^{-1}=D_4^{-1}(D_1+D_3)D_4^{-1}$. For any  nonzero vector $x$ 
in $\Re^n$, 
\begin{eqnarray*}
x^TF^{-1}x &=& x^TD_4^{-1}(D_1+D_3)D_4^{-1}x\\ 
&=&
x^TD_4^{-1}D_1D_4^{-1}x+x^TD_4^{-1}D_3D_4^{-1}x
\end{eqnarray*} 
Here we claim that $D_4^{-1}D_3D_4^{-1}$ is skew symmetric and 
therefore $x^TD_4^{-1}D_3D_4^{-1}x=0$. To obtain this one can 
notice that since $D_2=D_4^2$ and $D_2$ is a symmetric matrix, 
then $D_4$ and $D_4^{-1}$ 
are also symmetric.\\ 
Hence,
\[(D_4^{-1}D_3D_4^{-1})^T=D_4^{-1}D_3^TD_4^{-1}=-D_4^{-1}D_3D_4^{-1}.\] 
Thus 
\[x^TF^{-1}x= x^TD_4^{-1}D_1D_4^{-1}x, \hspace{.25in} \mbox{~for any~} 
\   0\ne x\in\Re^n.  \] Using the fact that $D_1-D_2$ is 
nonnegative, we obtain  
\[x^TF^{-1}x \ge x^TD_4^{-1}D_2D_4^{-1}x=x^Tx, \hspace{.25in} \mbox{~for any~} 
\   0\ne x\in\Re^n.  \] 
 This implies that 
\[{\norm{x}_2}^2 \le x^TF^{-1}x \le \norm{x}_2. \norm{F^{-1}x}_2, \hspace{.25in} \mbox{~for any~} \ 0\ne 
x\in\Re^n,  \] that is, 
\begin{equation} \label{Fbound}
\norm{x}_2 \le \norm{F^{-1}x}_2,\hspace{.25in} \mbox{~for any~} \   0\ne x\in\Re^n. 
\end{equation} 
Obviously (\ref{Fbound}) is equivalent to 
\begin{equation} 
\norm{Fy}_2 \le \norm{y}_2,\hspace{.25in} 
\mbox{~for any~} \   0\ne y\in\Re^n. 
\end{equation} 
Since the last equation holds for any nonzero vector $y$, then $\norm{F}_2 \le 1$. \\ 
\hspace{2in}$\Box{}{}$\\ 

For the next step, we analyze the stability of the nonhomogenous 
problem over an arbitrary 
but finite time interval $[0,T]$. We later show that the stability of 
the homogeneous problem does not depend on $T$. 
Consider (\ref{eq:schemeEquation}) with ${f}\not\equiv{0}$. 

After some simple calculations, we get that $u_n$ satisfies
\begin{equation} \label{eq:mainRecursion} 
u_{n+1}=(I+kA+kB_n)^{-1}(I+kC)u_n+k(I+kA+kB_n)^{-1}f_{n+1}. 
\end{equation} 

We denote the range 
of the step index $n$, by $[0,N]$, where $kN=T$. To simplify 
the notation, we do not explicitly indicate that $N$ depends on $k$
and $T$. 

\begin{theorem} \label{t:bouHom} 
Let (\ref{eq:condEquation}) hold. Then 
the solution of (\ref{eq:mainRecursion}) satisfies the following bound:
\begin{eqnarray*}  
\norm{u_{n+1}}_E 
&\le& 
\norm{u_0}_E+\frac{k}{1+k \lambda _{min}(C)}\Sigma_{p=0}^{n}\norm{f_{p+1}}_E\\
&\le&
\norm{u_0}_E+\frac{T}{(1+k\lambda_{min}(C))} \max_{t \in [0,T]} \norm{f(t)}_E,\, \forall 0 \leq n \leq N-1.
\end{eqnarray*}
Here $T$ is the size of the integration interval. 

\end{theorem} 
 \noindent 
{\bf Proof:} To simplify notation, we take $N_k=(I+kC)^{\frac{1}{2}}$
and $M_k=(I+kA+kB_n)^{-1}(I+kC)$. Then the equation (\ref{eq:mainRecursion})
can be written as 
\[u_{n+1}=M_ku_n+k(I+kA+kB_n)^{-1}f_{n+1}.\]
Using the definition 2.1, we have
\[(N_ku_{n+1})^T(N_ku_{n+1})=(N_ku_{n+1})^TN_kM_ku_n+k(N_ku_{n+1})
^TN_k(I+kA+kB_n)^{-1} 
f_{n+1}.\] 
\noindent 
Algebraic manipulation and the Cauchy-Schwarz inequality yield
\begin{eqnarray*}
{\norm{N_ku_{n+1}}_2}^2 
&\le& \norm{N_ku_{n+1}}_2 .\norm{N_kM_kN_k^{-1}}_2 . \norm{N_ku_n}_2\\
& &
+k\norm{N_ku_{n+1}}_2.\norm{N_kM_kN_k^{-1}}_2.\norm{N_k^{-1}f_{n+1}}_2.
\end{eqnarray*}
Using Lemma \ref{l:mainLemma} with $D_2=N_k^2$ and $D_1+D_3=M_kN_k^{-2}$, 
we obtain that $\norm{N_kM_kN_k^{-1}}_2 \le 1$. Then the  previous inequality reduces to
\[\norm{N_ku_{n+1}}_2 \le \norm{N_ku_n}_2+k\norm{N_k^{-1}f_{n+1}}_2\]
This inequality can be simplified as follows:
\begin{eqnarray*}
\norm{N_ku_{n+1}}_2 
&\le&
\norm{N_ku_n}_2 + k\norm{N_k^{-2}N_kf_{n+1}}_2\\
&\le&
\norm{N_ku_n}_2 + k\norm{(I+kC)^{-1}}_2 \norm{N_kf_{n+1}}_2\\ 
&\le&
\norm{N_ku_n}_2 + \frac{k}{(1+ k \lambda_{min}(C))}\norm{N_kf_{n+1}}_2.
\end{eqnarray*}
\noindent 
Thus,
\[\norm{u_{n+1}}_E - \norm{u_n}_E \le 
\frac{k}{(1+ k \lambda_{min}(C))}\norm{f_{n+1}}_E,\] 
\noindent and since $(I+kC)^{-1}$ is a symmetric positive definite 
matrix, 
\[\norm{I+kC}_2=\max{\lambda(I+kC)^{-1}} = \frac{1}{(\min{\lambda(I+kC)})}.\]
By the spectral 
mapping theorem $\lambda(I+kC)=1+k\lambda(C)$. Therefore 
\[\norm{(I+kC)^{-1}}_2= \frac{1}{1+k \lambda_{min}(C)},\]
\noindent where $\lambda_{min}(C)$ is the minimum eigenvalue of 
matrix $C$. This implies
\[\norm{u_{n+1}}_E - \norm{u_n}_E \le 
\frac{k}{(1+k \lambda _{min}(C))}\norm{f_{n+1}}_E,\, 0 \leq n \leq N-1. \] 
\noindent 
Summing from $0$ to $n$ gives
\[\norm{u_{n+1}}_E - \norm{u_0}_E \le \frac{k}{(1+k \lambda _{min}(C))}
\Sigma_{p=0}^{n}\norm{f_{p+1}}_E,\, \forall 0 \leq n \leq N-1,\] 
\noindent 
that is,
\[\norm{u_{n+1}}_E \le \norm{u_0}_E + \frac{k}{(1+k \lambda _{min}(C))}
\Sigma_{p=0}^{n}\norm{f_{p+1}}_E,\, 0 \leq n \leq N-1,\] 
which is the claimed first result. The second result follows immediately. 
$\Box{}{}$

The local truncation error of the method (\ref{eq:schemeEquation}) is clearly 
$O(\Delta t)$. In the error estimate (which follows) we need 
a precise statement of this fact, which we now derive.\\ 
To simplify  
our notation, we use $u_n$ to denote $u(t_n)$, where $u(\cdot)$ is the exact 
solution of (\ref{eq:mainEquation}). We also use $u_n$ to denote an iterate 
of our numerical scheme, but the particular meaning of $u_n$ will become 
evident from the context.  
 
According to the definition of local truncation error ~\cite{A89}, 
\begin{eqnarray}
\tau_{n+1}
&=&
\frac{u(t_{n+1})-u(t_n)}{k}+Au(t_{n+1})+B(u(t_n))u(t_{n+1})-Cu(t_n)\nonumber\\
& & 
-[u'(t_{n+1})+Au(t_{n+1})+B(u(t_{n+1}))u(t_{n+1})-Cu(t_{n+1})]\\
&=&
\frac{u_{n+1}-u_n}{k}-u'_{n+1}-(B(u(t_{n+1}))-B(u(t_n)))
u(t_{n+1})+ C(u_{n+1}-u_n). \nonumber 
\end{eqnarray} 
\noindent 
 
Using the second-order integral form of the Taylor expansion around 
$t_{n+1}$, we obtain  
\[  
u_{n+1}-u_n-k u'_{n+1}=-\int_{t_{n+1}}^{t_n} u''(t) (t-t_{n+1}) dt,  
\] 
which we rewrite as  
\[  
\frac{u_{n+1}-u_n}{k}- u'_{n+1}=-\frac{1}{k}\int_{t_{n+1}}^{t_n} u''(t) 
(t-t_{n+1}) dt= -\frac{1}{k}\int_{t_{n}}^{t_{n+1}} u''(t) (t_{n+1}-t) dt. 
\] 
Using the first-order integral form of the Taylor expansion around $t_n$, we  
obtain  
\[\begin{array}{rcl} 
& & (B(u(t_{n+1}))-B(u(t_n)))u(t_{n+1})-C(u_{n+1}-u_n) =  \\
& & \int_{t_n}^{t_{n+1}} \left( \frac{d}{dt} B(u(t)) u(t_{n+1})-C u'(t) 
\right) dt. 
\end{array}
\] 
 
Using the expression we have derived for the local truncation error 
$\tau_{n+1}$, and the preceding equations derived from Taylor's theorem, 
we obtain  
\begin{eqnarray*} 
\tau_{n+1} & = & -\frac{1}{k}\int_{t_{n}}^{t_{n+1}} u''(t) (t_{n+1}-t) dt -  
\int_{t_n}^{t_{n+1}} \left( \frac{d}{dt} B(u(t)) u(t_{n+1})-C u'(t) \right) dt 
\\ 
& = & \int_{t_n}^{t_{n+1}} \left(-\frac{t_{n+1}-t}{k} u''(t)- 
\frac{d}{dt} B(u(t)) u(t_{n+1})+C u'(t) \right) dt.  
\end{eqnarray*} 
  
By the mean value theorem, there exists $\xi_n \in (t_n,t_{n+1})$  
such that  
\begin{equation} \label{eq:mainTruncation} 
\tau_{n+1}  =  - u''(\xi_n) (t_{n+1}-\xi_n) - k \left.\frac{d}{dt} B(u(t)) 
\right|_{t=\xi_n} u(t_{n+1}) + k C u'(\xi_n).  
\end{equation} 
Hence, using the fact that $ 0 \leq \left(t_{n+1}-\xi_n \right) \leq k$, we 
obtain that  
\[ 
\norm{\tau_{n+1}}_2 \le k  
\max_{t_n \le s \le t_{n+1}}  
\left( 
\norm{ u''(s)}_2 + 
\norm{\left.\frac{d}{dt} B(u(t)) \right|_{t=s}}_2 
\max_{t_n \le \theta \le t_{n+1}}\norm{u(\theta)}_2+ 
\norm{Cu'(s)}_2 \right). 
\] 
\noindent 
This proves the following lemma. 
\begin{lemma} 
Let $k=\Delta t$ and $n\ge 0$. The method 
\begin{equation} \label{inhomogeneous}
\frac{u_{n+1}-u_n}{k}+Au_{n+1}+B_nu_{n+1}-Cu_n=f_{n+1},
\end{equation} 
\noindent where $A=A^T \succ 0$ and $C=C^T\succeq 0$ are ${n}\times{n}$ 
symmetric matrices, $B_n$ an ${n}\times{n}$ skew-symmetric matrix, 
and $f_{n+1}=f((n+1)k)$, is consistent. That is,  the local truncation 
error is $O(\Delta t)$. 
\end{lemma} 
 
We now bound the total error. We consider first the energy norm of truncation 
error. 
\begin{lemma} \label{l:btruncation}
Let $\tau_{n+1}$ be the local truncation error of method 
(\ref{inhomogeneous}). 
Then
\begin{equation}\label{eq:btruncation}
\norm{\tau_{n+1}}_E \le k  
\max_{0\le t \le T} 
\left(\norm{u''(t)}_E+ 
\norm{Cu'(t)}_E +  
\norm{\frac{d}{dt} B(u(t)) }_E 
\max_{0 \le s \le T} \norm{u(s)}_E \right).
\end{equation} 
\end{lemma} 
\noindent 
{\bf Proof:} By definition of energy norm and following the identity  
(\ref{eq:mainTruncation}), we get
\[\norm{\tau_{n+1}}_E =  
\norm{- u''(\xi_n) (t_{n+1}-\xi_n) -  
k \left.\frac{d}{dt} B(u(t)) \right|_{t=\xi_n}  
u(t_{n+1}) + k C u'(\xi_n)}_E \] 
for some $\xi_n \in [t_n, t_{n+1}]$.  
The conclusion follows after applying the inequality $0 \leq t_{n+1}-\xi_n 
\leq k$,  
the triangle inequality, and the properties of the $\max$ function.  
Note that $\norm{\frac{d}{dt} B(u(t)) }_E$ is the  
induced $\norm{\cdot}_E$ of the corresponding matrix.  
$\Box{}{}$ \\ 

We now give a convergence result for the solution of (\ref{eq:schemeEquation}). 
First, we need to compute a certain estimate. We have that 
\begin{eqnarray*}
& & \left[B(u(t_n))-B(u_n)\right] u(t_{n+1}) 
= 
\int_{0}^{1} \frac{d}{d \theta}[B(u(t_n) \theta + u_n(1-\theta))] u(t_{n+1}) dt = \\
&&
\left(\nabla_u (B(u(t_n) \theta_n + u_n(1-\theta_n)))e_n\right)u(t_{n+1}),\,
\mbox{~for some~} \theta_n \in [0,1],\\ 
\end{eqnarray*}
where $e_n=u(t_n)-u_n$. Here $u(t_n)$ is the 
solution of (\ref{eq:mainEquation}), whereas $u_n$ is the solution of 
our numerical scheme. 

We define the matrix $W_n$, by its action on a vector $x \in \Re^n$:
\[ W_n x = \left[ \left(\nabla_u B(u(t_n) \theta_n + u_n(1-\theta_n)) \right) x \right] u(t_{n+1}), \]
which results in the following identity
\begin{equation} \label{eq:wn}
 \left[B(u(t_n))-B(u_n)\right] u(t_{n+1}) = W_n e_n.
\end{equation}

\begin{lemma}\label{l:boundW}
Let u(.) be the solution of (\ref{eq:mainEquation}) and $u_n$ be
the approximation to $u(n\Delta t)$, obtained from the numerical scheme
 (\ref{eq:schemeEquation}). Then there exists $\Gamma$ such that, 
$\forall t \in [0,T]$ we have that 
\[ \norm{W_n}_2 \leq \Gamma,\mbox{~and~} \norm{W_n}_E \leq \Gamma_E=\Gamma \sqrt{1+k\lambda_{max}(C)},\;
\forall 0 \leq n \leq N.\]
\end{lemma}

{\bf Proof:} From Theorem \ref{t:bouHom} we have that 
\[ \begin{array}{rcl}
\norm{u_n}_2 & \leq & \norm{u_n}_E \leq \norm{u_0}_E + T \max_{t \in [0,T]} \norm{f(t)}_E \\
& \leq & \sqrt{1 + k \lambda_{max}(C)} \left( \norm{u_0}_2 + T \max_{t \in [0,T]} \norm{f(t)}_2 \right),\;
\forall 0\leq n \leq N. 
\end{array}
\]
We define 
\[ \Lambda_E=\sqrt{1+ T \lambda_{max}(C)} \left( \norm{u_0}_2 + T \max_{t \in [0,T]} \norm{f(t)}_2 \right).
\]
From Lemma \ref{l:boundedSolution} we have that $u(t)$ is bounded on 
$[0,T]$, and we define $\Lambda_u=\max_{t \in [0,T]} \norm{u(t)}_2$. Since 
$B(\cdot)$ is of class $C^1$, we can define 
\[
\Gamma=\max_{\theta \in [0,1],\, \norm{u_1}_2 \leq \Lambda_E,\, \norm{u_2}_2 =1,\, \norm{v_1} \leq \Lambda_u,
\norm{v}_2 \leq \Lambda_u} \norm{\left[\left(\nabla_u B(\theta v_2 + 
(1-\theta) u_1) \right) u_2\right] v_1}_2.
\]
From the definition of $W_n$, we immediately obtain that 
\[ 
\norm{W_n}_2 \leq \Gamma, \; \forall 0 \leq n \leq N.
\]
The second part of the conclusion follows from the inequality between 
$\norm{\cdot}_E$ and $\norm{\cdot}_2$. $\Box$

\begin{theorem}\label{t:main}
Consider solving the nonhomogenous problem on the interval [0,T]\\
\[u' +Au+B(u)u-Cu=f\]
\noindent
using the following method\\
\[\frac{u_{n+1}-u_n}{k}+Au_{n+1}+B_nu_{n+1}-Cu_n=f_{n+1},\]\\
\noindent where $k=\Delta t$, $B_n=B(u_n)$ and
$f_{n+1}=f((n+1)k)$. Let $e_n=u(t_n)-u_n$ denote the local error. 
Assume that $e_0=0$. Then the method is convergent and
\begin{eqnarray*}
\norm{e_{n+1}}_E 
& \le & 
\frac{\left(1+\frac{k \Gamma_E}{1+k \lambda _{min}(C)}\right)^{n-1}-1}
{1+\frac{k \Gamma_E }{1+k \lambda_{min}(C)}}\frac{k^2U}{1+k \lambda_{min}(C)} \\ 
& \le & 
\frac{e^\frac{ T \Gamma_E}{1+k \lambda _{min}(C)}-1}
{1+\frac{k \Gamma_E }{1+k \lambda_{min}(C)}}\frac{k^2U}{1+k \lambda_{min}(C)},\;
\forall 0 \leq n \leq N-1,
\end{eqnarray*}
when $\Gamma_E \neq 0$, and
\begin{eqnarray*}
\norm{e_{n+1}}_E 
\le 
(n+1) \frac{k^2U}{1+k \lambda_{min}(C)} 
\le 
T \frac{kU}{1+k \lambda_{min}(C)},\, 
\forall 0 \leq n \leq N-1,
\end{eqnarray*}
when $\Gamma_E=0$, where
\[ U = \max_{0\le t \le T}\left(\norm{u''(t)}_E
+\norm{Cu'(t)}_E + \norm{\frac{d}{dt} B(u(t)) }_E
\max_{0 \le s \le T} \norm{u(s)}_E  \right). \]
\end{theorem}
\noindent 
{\bf Proof:} Following the definition of the truncation error $\tau_{n+1}$ and 
using the equation (\ref{eq:wn}), 
we obtain that the error, $e_n=u(t_n)-u_n$, satisfies
\[\frac{e_{n+1}-e_n}{k}+Ae_{n+1}+B_ne_{n+1}-Ce_n=\tau_{n+1} - W_n e_n.\] 
\noindent 
After algebraic calculations, we find that
\[e_{n+1}=(I+kA+kB_n)^{-1}(I+kC)e_n+k(I+kA+kB_n)^{-1}(\tau_{n+1}-W_n e_n).\] 
\noindent 
We use the energy inner product to obtain
\[\begin{array}{l} <e_{n+1},e_{n+1}>_E= \\ 
<(I+kA+kB_n)^{-1}(I+kC)e_n+k(I+kA+kB_n)^{-1}
(\tau_{n+1}-W_n e_n),e_{n+1}>_E.
\end{array} \] 
\noindent 
Applying the definition of energy norm (\ref{def:ENorm}) and  
the substitutions $M_k=(I+kA+kB_n)^{-1}(I+kC)$, and $N_k=(I+kC)^\frac{1}{2}$, we find that  
\[\begin{array}{l} (N_ke_{n+1})^T(N_ke_{n+1})= \\
(N_ke_{n+1})^TN_kM_ke_n+k(N_ke_{n+1})^TN_k
(I+kA+kB_n)^{-1})(\tau_{n+1}-W_n e_n). 
\end{array}
\] 
\noindent 
Using the Cauchy-Schwarz inequality, we obtain that 
\begin{eqnarray*}
\norm{N_ke_{n+1}}_2^2 
&\le&
\norm{N_ke_{n+1}}_2 .\norm{N_kM_kN_k^{-1}}_2 . \norm{N_ke_n}_2\\
& &
+k\norm{N_ke_{n+1}}_2.\norm{N_kM_kN_k^{-1}}_2 .\norm{N_k^{-1} 
\left(\tau_{n+1}-W_n e_n \right)}_2.
\end{eqnarray*} 
\noindent 
Thus 
\[\norm{N_ke_{n+1}}_2 \le \norm{N_kM_kN_k^{-1}}_2 . \norm{N_ke_n}_2+ 
k\norm{N_kM_kN_k^{-1}}_2 .\norm{N_k^{-1} \left(\tau_{n+1}-W_n e_n \right)}_2.\] 
\noindent 
Using Lemma \ref{l:mainLemma} with $D_2=N_k^2$ and $D_1+D_3=M_kN_k^{-2}$, 
we obtain that $\norm{N_kM_kN_k^{-1}}_2 \le 1$. 
\noindent 
Hence
\[\norm{e_{n+1}}_E \le \norm{e_n}_E + 
k\norm{N_k^{-2}}_2 \norm{ \left(\tau_{n+1}-W_n e_n \right)}_E.\] 
\noindent 
Equivalently, we obtain that 
\[\norm{e_{n+1}}_E \le \norm{e_n}_E + 
k\norm{(I+kC)^{-1}}_2 \left(\norm{\tau_{n+1}}_E + \norm{W_n}_E 
\norm{e_n}_E\right).\]
\noindent Notice that $(I+kC)^{-1}$ is a symmetric 
positive definite matrix and 
\[\norm{(I+kC)^{-1}}_2=\frac{1}{1+k \lambda _{min}(C)}.\]
On the other hand, by Lemma \ref{l:boundW}, there is a constant
$\Gamma_E$ such that $\norm{W_n}_E \le \Gamma_E$.
Therefore,
\begin{equation}\label{eq:recursion}
\norm{e_{n+1}}_E \le \left(1+\frac{k\Gamma_E}{1+k \lambda _{min}(C)}
\right)\norm{e_n}_E + \frac{k}{1+k \lambda _{min}(C)}\norm{\tau_{n+1}}_E.
\end{equation} 

This is a recursion formula of the following form: 
\[ 
r_{n+1} \leq a r_n + b \tau_n, \\
\]
which, when $a \neq 0$  has an upper bound of the type 
\[ r_{n+1} \leq a^n r_0 + \frac{a^{n-1}-1}{a} b \max_n \norm{\tau_{n}}_E. \]

Using this fact, we obtain that, when $\Gamma_E \neq 0$, the following bound for the error holds 
whenever $0 \leq n \leq N-1.$

\begin{eqnarray*}
\norm{e_{n+1}}_E 
&\le&
\left(1+\frac{k\Gamma_E}{1+k \lambda _{min}(C)}\right)^n\norm{e_0}_E\\
& &
+\frac{\left(1+\frac{k\Gamma_E}{1+k \lambda _{min}(C)}\right)^{n-1}-1}
{1+\frac{k\Gamma_E}{1+k \lambda_{min}(C)}}.\frac{k}{1+k \lambda _{min}(C)}
\max_n \norm{\tau_{n+1}}_E
\end{eqnarray*}

Replacing $\norm{\tau_{n+1}}_E$ by its bound 
(\ref{eq:btruncation}) obtained in Lemma \ref{l:btruncation}, 
and considering that $e_0=0$, we have, when $\Gamma_E \neq 0$ and $0 \leq n \leq N-1$, that

\begin{equation*}
\norm{e_{n+1}}_E 
\le
\frac{\left(1+\frac{k\Gamma_E}{1+k \lambda _{min}(C)}\right)^{n-1}-1}
{1+\frac{k\Gamma_E}{1+k \lambda _{min}(C)}}.\frac{k^2U}{1+k \lambda _{min}(C)}
\end{equation*}
with $U=\max_{0\le t \le T}\left(\norm{u''(t)}_E
+\norm{Cu'(t)}_E + \norm{\frac{d}{dt} B(u(t)) }_E
\max_{0 \le s \le T} \norm{u(s)}_E  \right)$.  
The second inequality for $\Gamma \neq 0$ follows from the inequality 
$(1+x)^n \leq e^{xn}$, for $x>0$ and $n$ positive integer. 

When $\Gamma_E = 0$, we immediately get from (\ref{eq:recursion}) and 
from Lemma \ref{l:btruncation} that 
\[
\norm{e_{n+1}}_E 
\le 
(n+1) \frac{k^2U}{1+k \lambda_{min}(C)},\; \forall 0 \leq  n \leq N-1, \]
which, together with $k N = T$ prove the inequalities for $\Gamma_E = 0$. 

The convergence 
follows from the fact that $\norm{\cdot}_E$ converges
to $\norm{\cdot}_2$ as $k \rightarrow 0$ which implies 
that $\norm{e_n}_2 \rightarrow 0$ as $k \rightarrow 0$. 
$\Box{}{}$
 
The case $\Gamma_E=0$ occurs, for example, when $B(u)$ is constant (which
we simulate numerically in the next section). For that case, the error 
increases only linearily with the size of the interval, assuming that 
the derivatives up to order $2$ of the solution $u(t)$ are uniformly bounded.  

\section{Numerical Results}

Let $\Omega=[0,1] \times [0,1]$. For the
equation 
\begin{eqnarray} \label{eq:modelPDE} 
u_t + b\cdot \nabla u - \epsilon \Delta u &=& f, \mbox{~over~}
\Omega,\nonumber
\\
u &=& \phi(x) \mbox{~on~} \delta \Omega, \\
u(x,0) &=& u_0(x) \mbox{~in~} \Omega, \nonumber
\end{eqnarray}
use the method described in this work, with uniform mesh and central 
difference. A choice must be made for the antidiffusion operator: 
averaging or projection. We have selected averaging. Since it is just
outside the theory, we will thereby test the robustness of the algorithm.    
Antidiffusion is completed by averaging, where
$\bar{u}(p)$:=weighted average of nearest neighbors. This
corresponds to filtering with $\delta=2h$. The method becomes in
our case
\[\dot u_{ij}(t) + b \cdot \nabla^h u_{ij} - (\epsilon + \epsilon_0) 
\Delta^h u_{ij} +
\epsilon_0 \overline{\Delta \overline{u_{ij}}}^q=h,\] where $q$ denotes 
how many times
the average operation is taken. In our experiments we chose $q=2$ and 
$\epsilon=10^{-4}$. We take
$b=(\cos(\theta),\sin(\theta))$, where $\theta=17^\circ$. 

For the boundary and initial conditions we take the line at angle
$\theta$ through the center of the domain. On the north side of
the line we take $\phi=1$ on the boundary; on the south side we
take $\phi=0$ on the boundary. We take $f=0$ and $0$ as initial
conditions. 

\begin{figure} 
\centerline{\psfig{figure=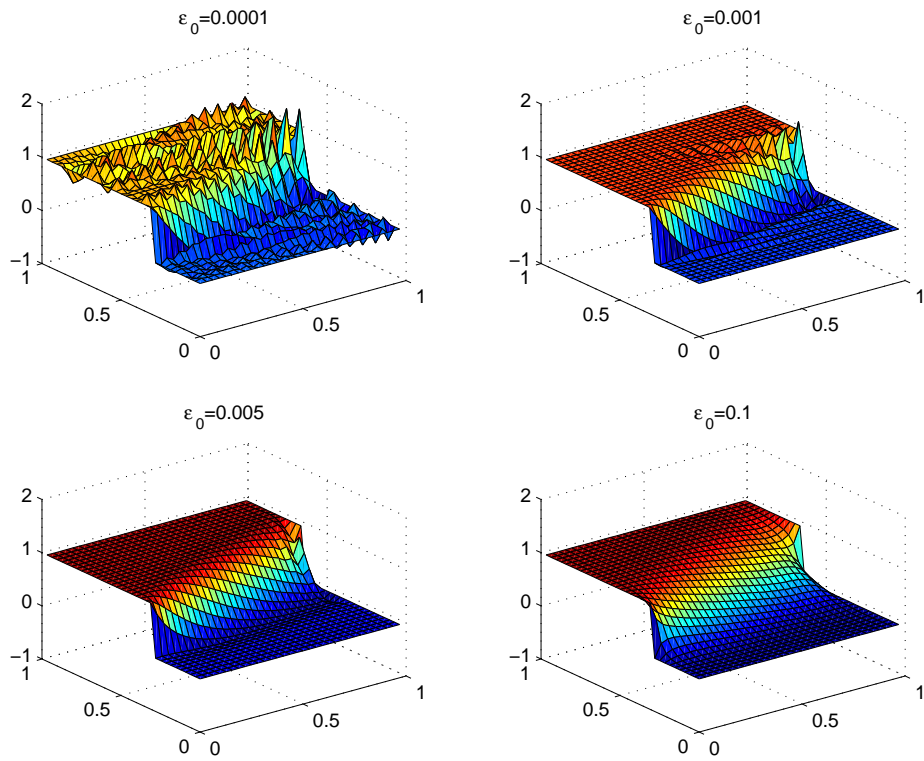,height=4in}}
\caption{\label{fig:surfaces}
Spatial stability of the steady-state solution for various choices of the 
artificial viscosity parameter $\epsilon_0$. }
\end{figure}

\begin{figure} 
\centerline{\psfig{figure=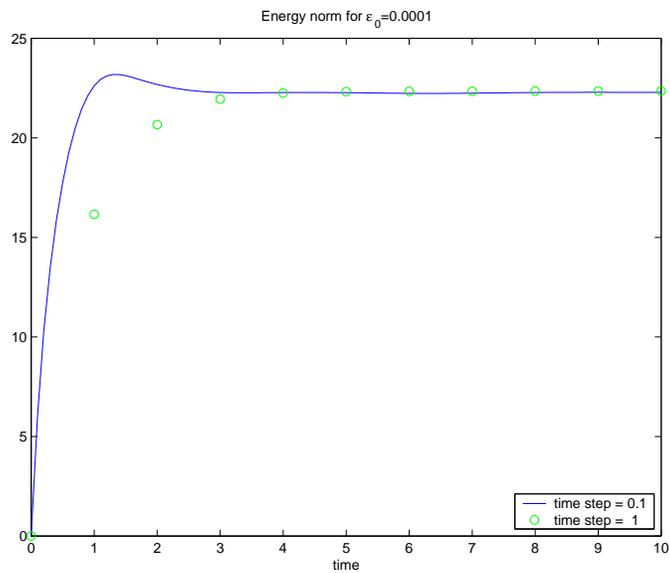,height=3in}}
\caption{\label{fig:timestepEnergyComparison} 
Stability of the numerical method demonstrated by the behavior of the  
energy norm}
\end{figure}

\begin{figure} 
\centerline{\psfig{figure=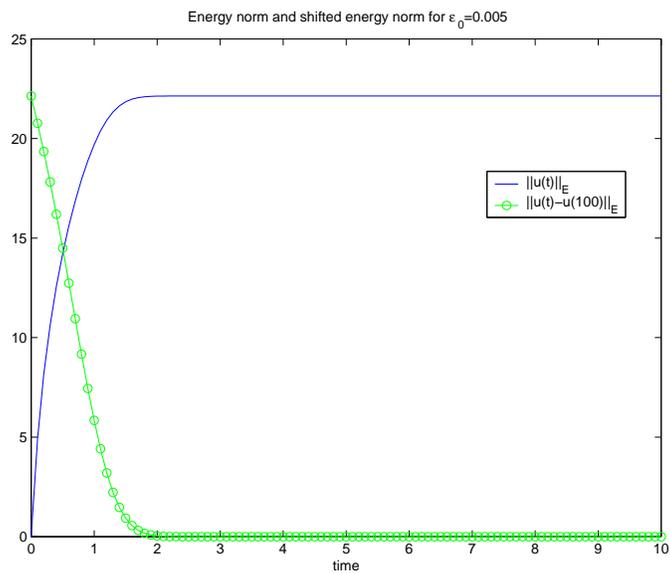,height=3in}}
\caption{\label{fig:shiftedBest} 
Numerical validation of Theorem \ref{t:unconditional}}
\end{figure}

\begin{figure} 
\centerline{\psfig{figure=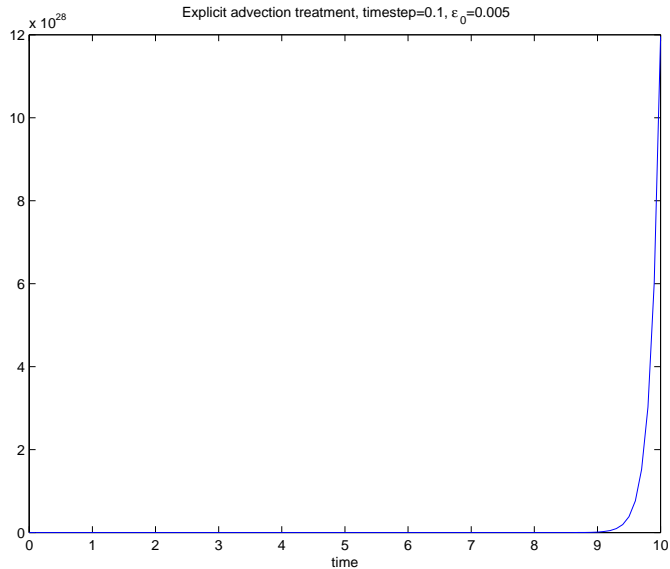,height=3in}}
\caption{\label{fig:explosion} 
Exponential growth of the solution of the scheme that includes the advection term 
explicitly}
\end{figure}

We performed the following experiments, all on a 32 $\times$ 32 mesh.
\begin{enumerate}
\item We ran the simulation for $1,000$ steps with a timestep of $10$, with the artificial viscosity 
parameter $\epsilon_0$ having succesively the values 
$10^{-1}$, $5 \times 10^{-3}$,
$10^{-3}$, $10^{-4}$. We have presented no analysis for the spatial dependence of the solution with respect to $\epsilon_0$, but we have 
included this experiment for validation, since our choice of parameters should result roughly in the steady-state 
approximation for this mesh, which has been studied before in the literature. 

The results are depicted in Figure \ref{fig:surfaces}. We see that when the artificial viscosity parameter 
$\epsilon_0$ is very small, a complete loss of coherence of the spatial structure results, 
whereas too large a parameter ($\epsilon_0=0.1$) 
alters the steady-state solution significantly. This effect is consistent 
with the typical behavior of centered methods for the skew step problem 
~\cite{RST96}. 
\item For $\epsilon_0=10^{-4}$, we ran the simulation for $100$ steps with a timestep of $1$ and for 
$1,000$ steps 
with a timestep of $0.1$. The energy norm comparison of these computations is 
presented in Figure \ref{fig:timestepEnergyComparison}.
We see that even for the very large step, the energy norm stays bounded, consistent with our absolute stability claim. 

We also present in Figure \ref{fig:shiftedBest} a comparison between
the energy norms of the distance between the successive iterates of
the two cases and their outcome at time $100$. From Figure
\ref{fig:timestepEnergyComparison}, we infer that $u(100)$ is a
reasonable approximation to the steady-state solution. Since the
equation (\ref{eq:modelPDE}) is linear, we have that $u_n-u(100)$ is
the result of the numerical scheme applied to the homogeneous equation
associated to (\ref{eq:modelPDE}). From Theorem \ref{t:unconditional}
we have that $\norm{u_n-u(100)}_E$ must be a decreasing sequence,
which is exactly what we observe from Figure
\ref{fig:shiftedBest}. Note that $\norm{u_n}_E$ is not a decreasing
sequence, as can be seen in Figure \ref{fig:shiftedBest}. Moreover,
the sequence $\norm{u_n}_E$ may not even be monotonic, as seen in
Figure \ref{fig:timestepEnergyComparison}, for $k=0.1$.

\item We compare the results of our scheme with the similar scheme that takes into account explicitly the 
term that contains the skew-symmetric matrix $B(u_n)$. For the latter scheme we obtain the recursion 
\[ 
\frac{u_{n+1}-u_n}{k}+Au_{n+1}+B(u_n)u_{n}-Cu_n=f_{n+1}. 
\]
We apply this scheme to our example on a 32 $\times$ 32 mesh for $1000$ timesteps of length $k=1$. 
We see the rapid exponential growth that is typical for computations with the timestep outside the region 
of stability. 

This demonstrates that our scheme has significantly better stability properties than the alternative, which 
would result in linear systems of comparable sparsity. The numerical scheme, based on a backward Euler approach that considers all terms implicitly, 
though absolutely stable, will result in less sparse linear systems since the matrix $C$ contains an averaging 
operator that substantially reduces sparsity and is not considered here for comparison.

\end{enumerate}

\section*{Acknowledgements}
This research was supported by the Department of 
Energy, through the Contract W-31-109-ENG-38, (MA), 
and the National Science Foundation through awards DMS-0112239 (MA and WJL).
and DMS-0207627 (FP and WJL).

\vfill
\begin{flushright}
\scriptsize
\framebox{\parbox{2.4in}{The submitted manuscript has been created
by the University of Chicago as Operator of Argonne
National Laboratory ("Argonne") under Contract No.\
W-31-109-ENG-38 with the U.S. Department of Energy.
The U.S. Government retains for itself, and others
acting on its behalf, a paid-up, nonexclusive, irrevocable
worldwide license in said article to reproduce,
prepare derivative works, distribute copies to the
public, and perform publicly and display publicly, by or on
behalf of the Government.}}
\end{flushright}

\bibliographystyle{siam} 
\bibliography{Xbib} 

\begin{thebibliography}{10}

\bibitem{A89}
{\sc K.~E. Atkinson}, {\em An introduction to numerical analysis}, Wiley, 1989.

\bibitem{BR62}
{\sc G.~Birkhoff and G.-C. Rota}, {\em Ordinary Differential Equations}, Ginn
  and Company, Boston, 1962.

\bibitem{G99a}
{\sc J.~L. Guermond}, {\em Stabilization of {G}alerkin approximations of
  transport equations by subgrid modeling}, M2AN, 33 (1999), pp.~1293--1316.

\bibitem{G99b}
\leavevmode\vrule height 2pt depth -1.6pt width 23pt, {\em Stabilization par
  viscosite de sous-maille pour l'approximation de {G}alerkin des operateurs
  lineaires monotones}, C.R.A.S., 328 (1999), pp.~617--622.

\bibitem{HMJ00}
{\sc T.~J. Hughes, L.~Mazzei, and K.~E. Jasen}, {\em Large eddy simulation and
  the variational multiscale method}, Comput.Visual Sci., 3 (2000), pp.~47--59.

\bibitem{HU00}
{\sc T.~J.~R. Hughes, L.~Mazzei, and K.~E. Jansen}, {\em Large eddy simulation
  and the variational multiscale method}, Comput. Visual Sci., 3 (2000),
  pp.~47--59.

\bibitem{IL98}
{\sc T.~Iliescu and W.~Layton}, {\em Approximating the largger eddies in fluid
  motion {III}: the {B}oussinesq model for turbulent fluctuations}, Analele
  Stiintifice ale Universitatii Al.l.Cuza, tomul XLIV (1998), pp.~245--261.

\bibitem{KA02}
{\sc S.~Kaya}, {\em Numerical analysis of a subgrid scale eddy viscosity method
  for higher reynolds number flow problem}, University of Pittsburgh,Technical
  report,  (2002).

\bibitem{KL02}
{\sc S.~Kaya and W.~Layton}, {\em Subgrid-scale eddy viscosity methods are
  variational multiscale methods}, University of Pittsburgh,Technical report,
  (2002).

\bibitem{KP80}
{\sc H.~Kesten and G.~Papanicolaou}, {\em A limit theorem for stochastic
  acceleration}, Comm. Math. Phys., 78 (1980), pp.~19--63.

\bibitem{L00}
{\sc W.~Layton}, {\em Approximating the larger eddies in fluid motion {V}:
  {K}inetic energy balance of scale similarity models}, Math. and Computer
  Modeling, 31 (2000), pp.~1--7.

\bibitem{L02}
\leavevmode\vrule height 2pt depth -1.6pt width 23pt, {\em A connection between
  subgrid scale eddy viscosity and mixed methods}, Appl. Math. and Computing,
  133 (2002), pp.~147--157.

\bibitem{TA89}
{\sc Y.~Maday and E.~Tadmor}, {\em Analysis of the spectral vanishing viscosity
  method for periodic conservation laws}, SIAM Journal on Numerical Analysis,
  26 (1989), pp.~854--870.

\bibitem{MP93}
{\sc B.~Mohammadi and O.~Pironneau}, {\em Analysis of the K-$\epsilon$
  Turbulence Model}, Wiley, 1993.

\bibitem{RST96}
{\sc H.~G. Roos, M.~Stynes, and L.~Tobiska}, {\em Numerical Methods for
  Singularly Perturbed Differential Equations}, Springer, Berlin, 1996.

\end{thebibliography}
\end{document}